\documentclass{article}
 \usepackage{amsmath}

\newtheorem{thm}{Theorem}[section]
\newtheorem{cor}[thm]{Corollary}
\newtheorem{prop}[thm]{Proposition}

\newtheorem{lem}[thm]{Lemma}

\newtheorem{rem}[thm]{Remark}

\newcommand{\be}{\begin{equation}}
\newcommand{\ee}{\end{equation}}
\newcommand{\ben}{\begin{enumerate}}
\newcommand{\een}{\end{enumerate}}
\newcommand{\beq}{\begin{eqnarray}}
\newcommand{\eeq}{\end{eqnarray}}
\newcommand{\beqn}{\begin{eqnarray*}}
\newcommand{\eeqn}{\end{eqnarray*}}

\newcommand{\pa}{\partial}

\newcommand{\qed}{\hspace*{\fill}Q.E.D.}  

\begin{document}
\title{On a Class of Complete and Projectively Flat Finsler Metrics}
\author{Guojun Yang\footnote{Supported by the National Natural Science Foundation of China (11471226)}}
\date{}
\maketitle
\begin{abstract}

An $(\alpha,\beta)$-manifold $(M,F)$ is a Finsler manifold with
the Finsler metric $F$ being  defined by a Riemannian metric
$\alpha$ and $1$-form $\beta$ on the manifold $M$. In this paper,
we classify $n$-dimensional $(\alpha,\beta)$-manifolds
(non-Randers type) which are positively complete and locally
projectively flat. We show that the non-trivial class is that $M$
is homeomorphic to the $n$-sphere $S^n$ and $(S^n,F)$ is
projectively related to a standard spherical Riemannian  manifold,
and then we obtain some special geometric properties on the
geodesics and scalar flag curvature of $F$ on $S^n$, especially
when $F$ is a metric of general square type.


\end{abstract}

\section{Introduction}

In Finsler geometry, the flag curvature is a natural extension of
the sectional curvature in Riemannian geometry.  For a Finsler
metric on an manifold $M$, the flag curvature $K= K(P,y)$ usually
depends on a section  (flag) $P\subset T_xM$ and a vector (pole)
$y\in P$. It  is said to be of scalar flag curvature if $K=K(x,y)$
is a scalar function of
 $y\in T_xM$, and of constant flag curvature if $K=constant$.  In dimension
  $n\geq 3$, a Riemannian metric is of scalar
  flag curvature if and only if it is of constant sectional curvature. For
 a Finsler manifold $(M,F)$, geodesics are locally minimizing
 curves on the manifold $M$ as to the distance defined by the
 metric $F$. Geodesics are the solutions of an ODE of second
 order. For a point $x\in M$ and a vector $y\in T_xM$, there is a
 unique geodesic $\gamma=\gamma(t)$ satisfying
 $\gamma(0)=x,\gamma'(0)=y$, where $t$ is the arc-length
 parameter of $\gamma$. Flag curvatures and geodesics are closely
 related through the second-variation of arc-length for geodesics,
 and both of them play an important role in the studies of Finsler
 geometry.

It is the Hilbert's Fourth Problem to study  locally projectively
flat metrics. A Finsler metric is called locally projectively flat
if its geodesics are locally straight. The  Beltrami Theorem
states that a Riemannian metric is locally projectively flat if
and only if it is of constant sectional curvature. For locally
projectively flat  Finsler metrics, they are generally not of
constant flag curvature, but they are always of scalar flag
curvature.  In \cite{Shen2}, Z. Shen classifies locally
projectively flat Finsler metrics of constant flag curvature by
using the model of Funk metric. A Funk metric, defined on a
strongly convex domain $\Omega\in\mathbf{R}^n$, is projectively
flat with negative constant flag curvature, and particularly, it
is positively complete (\cite{Ok}).  When $\Omega$ is the unit
ball $B^n$, the Funk metric becomes a Randers metric. In
\cite{SY}, it gives a class of square metrics which are positively
complete, locally projectively flat and of constant flag curvature
on the unit ball $B^n$. In this paper, we will classify a class of
(positively) complete Finsler metrics which are locally
projectively flat and show  some geometrical properties of their
geodesics and scalar flag curvatures.

All Finsler metrics on a manifold $M$ in this paper are assumed to
be regular, namely, they are positively definite and defined on
the whole slit tangent bundle $TM-0$. An $(\alpha,\beta)$-metric
$F$ is a Finsler metric  defined by a Riemannian metric
 $\alpha=\sqrt{a_{ij}(x)y^iy^j}$ and a $1$-form $\beta=b_i(x)y^i$ on a manifold
 $M$, which can be expressed in the following form:
 $$F=\alpha \phi(s),\ \ s=\beta/\alpha,$$
where $\phi(s)$ is a function defined on an open interval
$(-b_{o},b_o)$ (see the regular condition of $F$ in Section
\ref{sec2}). For a pair of $\alpha$ and $\beta$, define
$b:=||\beta||_{\alpha}$. If $\phi(s)=1+s$, then $F=\alpha+\beta$
with $b<1$ is called a Randers metric, which has a lot of special
properties. Two $(\alpha,\beta)$-metrics $F$ and $\widetilde{F}$
are called of the same metric type if
 $$F=\alpha\phi(\beta/\alpha),\
 \widetilde{F}=\widetilde{\alpha}\phi(\widetilde{\beta}/\widetilde{\alpha}):\
 \ \
 \widetilde{\alpha}=\sqrt{\alpha^2+\epsilon \beta^2},\
\widetilde{\beta}=k\beta,$$
 where $\epsilon,k$ are constant. In this paper, we mainly consider the
 $(\alpha,\beta)$-metric $F=\alpha\phi(\beta/\alpha)$ with $\phi(s)$ being defined by
 \be\label{y1}
    \big\{1+(k_1+k_3)s^2+k_2s^4\big\}\phi''(s)=(k_1+k_2s^2)\big\{\phi(s)-s\phi'(s)\big\},\
    \ \ (\phi(0)=1,\ k_2\ne k_1k_3),
   \ee
where $k_1,k_2,k_3$ are constant.  If $k_2=k_1k_3$, then $F$ is of
Randers metric
 type.  A special metric type of (\ref{y1}), called a general square metric type,
 is defined by
  \be\label{y2}
 F=F^{\pm}_\epsilon:=\alpha+\epsilon\beta\pm \frac{\beta^2}{\alpha},\ \ \
 (\epsilon=constant),
  \ee
 in which $F_2^+$ is called a square metric in some literatures (\cite{LS} \cite{SYa} \cite{SY}),
 and the metric $F_0^{\pm}$ also has some special properties (\cite{Y1}). Note that $F$ in (\ref{y1}) is of the metric
 type $F_0^{\pm}$ if and only if $k_2=(2k_1+3k_3)(3k_1+2k_3)/25$  (\cite{Y1}). It
 is shown in \cite{Shen1} that, in dimension $n\ge 3$, if an
 $(\alpha,\beta)$-metric $F$ (non Randers type, $\beta$ not
 parallel) is locally projectively flat, then $\phi(s)$ satisfies
 (\ref{y1}) and $\beta$ must be closed; while  in dimension $n=2$, there is one more class with
 $F$ being of the metric type $F_0^{\pm}$ and $\beta$ generally being not closed
 (\cite{Y1}). We  show that in dimension $n\ge 3$, the metric
 $F$ in (\ref{y2}) is locally projectively flat iff. $F$
 is of scalar flag curvature (\cite {SYa} \cite{Y3}).

\begin{thm}\label{th1}
Let $(M,F)$ be an  $n(\ge 2)$-dimensional Finsler manifold, where
$F=\alpha\phi(\beta/\alpha)$ is an $(\alpha,\beta)$-metric on $M$.
Assume that (C1) $F$ is of non-Randers metric type; (C2) $n\ge 3$
if $F$ is of the metric type $F_0^{\pm}$ in (\ref{y2}); (C3)
$(M,F)$ is positively complete; (C4) $F$ is locally projectively
flat; (C5) $b(x_0)=Sup_{x\in M} b(x)$ at some point $x_0\in M$.
Then one of the following cases holds:

 \ben
       \item[{\rm (i)}]  $F=\alpha$ is a complete Riemannian metric of constant sectional curvature.

       \item[{\rm (ii)}]   $F$ is  flat-parallel ($\alpha$ is flat and $\beta$ is parallel
       with respect to $\alpha$).

      \item[{\rm (iii)}] $M$ is homeomorphic to the $n$-sphere
      $S^n$ in the $(n+1)$-dimensional Euclidean space, and  ($S^n,F$) is projectively related to a
      Riemannian  metric on $S^n$ of positive constant sectional
      curvature. So all geodesics of ($S^n,F$) are closed.
\een
 Further, if $F$ satisfies the conditions (C3)--(C5), and
 (C6) $F$ is of constant flag curvature, then only cases
(i) and (ii) occur.
\end{thm}

The condition (C5) can be replaced by a weaker condition (see
Remark \ref{rem31} below). On compact manifold (without boundary),
(C3) and (C5) are obviously satisfied. Without the condition (C5),
there are Randers metrics or square metrics which satisfy (C3),
(C4) and (C6) but do not belong to the cases (i), (ii) or (iii) of
Theorem \ref{th1} (\cite{Shen0} \cite{SY}). Cases (i) and (ii) of
Theorem \ref{th1} are trivial. So below we further study case
(iii) of Theorem \ref{th1}, and in this case, the metric $F$ is
determined by the ODE (\ref{y1}).

Corresponding to the ODE (\ref{y1}), let $u=u(t),v=v(t),w=w(t)$
satisfy the following ODEs:
 \beq
u'&=&\frac{v-k_1u}{1+(k_1+k_3)t+k_2t^2},\label{y3}\\
v'&=&\frac{u(k_2u-k_3v-2k_1v)+2v^2}{u[1+(k_1+k_3)t+k_2t^2]},\label{y4}\\
w'&=&\frac{w(3v-k_3u-2k_1u)}{2u\big[1+(k_1+k_3)t+k_2t^2\big]}.\label{y5}
 \eeq
Make a transformation
 \be\label{y6}
h:=\sqrt{u\alpha^2+v\beta^2},\ \ \ \rho:=w\beta,
 \ee
where $u=u(b^2)>0,v=v(b^2),w=w(b^2)\ne 0$ satisfy the ODEs
(\ref{y3})--(\ref{y5}) such that $h$  is a Riemann metric and the
1-form $\rho$ satisfies certain norm condition with respect to
$h$. There are different choices for $u,v,w$. It can be  shown
 that, if $F$ defined by (\ref{y1}) is locally
projectively flat, then $h$ is of constant sectional curvature
$\mu$ and $\rho$ is closed and conformal with respect to $h$ (cf.
\cite{Y3} \cite{Yu1}). Then the covariant
    derivatives $p_{i|j}$ of $\rho=p_iy^i$ with
    respect to $h=\sqrt{h_{ij}y^iy^j}$ satisfy
 \be\label{y7}
 p_{i|j}=-2ch_{ij}
 \ee
 for some scalar function $c=c(x)$. In the whole paper, we will
 use the data $(u,v,w)$ and $(h,\mu,c,\delta)$ (related to
  the choice of $(u,v,w)$), where  $\delta$ is
 defined in Theorem \ref{th2} (i) below.

\begin{thm}\label{th2}
 Let $(S^n,F)$ be an  $n(\ge 2)$-dimensional non-Riemannian Finsler manifold, where
$F=\alpha\phi(\beta/\alpha)$ is an $(\alpha,\beta)$-metric defined
by the ODE (\ref{y1}).
 Assume the
 conditions  (C2) and (C4) in Theorem \ref{th1} are satisfied.
 Then the following hold

 \ben
     \item[{\rm (i)}] For every suitable choice of $u,v$ and $w$, we
     have $\rho=2\mu^{-1}c_0$, where $c_0:=c_iy^i=c_{x^i}y^i$.
     Further, the gradient field $\nabla c$ has just two vanishing points $P,Q\in S^n$
     and  $c=\mu^{-\frac{1}{2}}\delta \cos(\sqrt{\mu}\ t)$,  where $t$ is the arc-length parameter of any geodesic
      of $(S^n,h)$ connecting $P$ and $Q$, and
     $\delta:=\sqrt{||\nabla c||_h^2+\mu c^2}$ is a constant.

    \item[{\rm (ii)}] For a special suitable choice of $u,v$ and $w$, the
    arc-length
    $L$ of any closed geodesic of $(S^n,F)$ through $P,Q$ and  any  closed geodesic
    of $(S^n,F)$ on the  hypersurface  $c=0$ has  the following  expansion   as to $\delta$
   \be\label{y8}
   L=\frac{2\pi}{\sqrt{\mu}}-\frac{4k_3\pi}{\mu^2\sqrt{\mu}}\ \delta^2+
   \frac{4(3k_3^2+2k_1k_3-2k_2)\pi}{\mu^4\sqrt{\mu}}\ \delta^4+o(\delta^6),\\
   \ee
 When $\delta\rightarrow 0$ is small enough, $L>2\pi/\sqrt{\mu}$
  if $k_3<0$ or $k_3=0$ but $k_2<0$; and $L<2\pi/\sqrt{\mu}$
  if $k_3>0$ or $k_3=0$ but $k_2>0$.

 \item[{\rm (iii)}] The maximal and minimal values of
 the scalar flag curvature $K$ of  $F$ on $S^n$ are
 just that of the function $R=R(s,t)$   of two variables $(s,t)$ defined on a
   bounded and  closed subset  $D=\{(s,t)|0\leq
   t\leq \delta^2/\mu,\  s^2\leq B\}$ in the Euclidean plane, where
   $R$ is defined by
  \be\label{y9}
 R=\frac{u}{\phi^2
 w^2}\Big\{3u\big(\frac{f_1^2\phi'^2}{\phi^2}+\frac{4f_1-f_2^2}{s^2}\big)
 +2(uf_3-v)\big(f_2-\frac{f_1s\phi'}{\phi}\big)\Big\}t+
 \frac{\mu u (f_1s\phi'-f_2\phi)}{2\phi^3}
  \ee
 for every suitable choice of $u,v,w$ (functions of $B$), in which $f_1,f_2,f_3$ are defined by
   $$
   f_1:=1+(k_1+k_3)s^2+k_2s^4, \ \ f_2:=k_2s^4-k_1s^2-2, \ \
 f_3:=3k_2s^2+k_1+3k_3,
   $$
   and $B=B(t):=b^2$ is the unique solution of  the following
   equation for $0\leq  t\leq \delta^2/\mu$
   \be\label{y10}
    \frac{w^2(B)B}{u(B)+v(B)B}=\frac{4(\delta^2-\mu
       t)}{\mu^2}.
   \ee
   Further, the scalar flag curvature $K$ is related to the function $R$ by $K=R(\beta/\alpha,c^2)$.
 \een
\end{thm}

In the proof of Theorem \ref{th2} (ii), we will give the integral
expression of (\ref{y8}). The expansion (\ref{y8}) is only for two
families of closed geodesics of $(S^n,F)$, and for other closed
geodesics of $F$ on $S^n$, we have not gotten the estimation of
their arc-length. For a general metric in Theorem \ref{th2} (iii),
it is difficult to obtain the explicit maximal and minimal values
of $R(s,t)$ on $D$. By a result in \cite{MS}, it is seen that
$Max_{(s,t)\in D}\big(R(s,t)\big)$ can never be negative. In
Remark \ref{rem42} below, we show a little different function
$\widetilde{R}(s,t)$ on $\widetilde{D}$, which is equivalent to
$R(s,t)$ on $D$. For some general square metrics in (\ref{y2}),
there is a simpler estimation on geodesics and the scalar flag
curvature $K$ by a different choice of $u,v$ and $w$.

\begin{thm}\label{th3}
  Let $(S^n,F)$ be an  $n(\ge 2)$-dimensional non-Riemannian Finsler manifold, where
$F=F_{\epsilon}^{\pm}$ is a general square metric defined by
 (\ref{y2}).
 Assume the
 conditions  (C2) and (C4) in Theorem \ref{th1} are satisfied.
 Let $u=(1\mp b^2)^2,\ v=0, \ w=\sqrt{1\mp  b^2}$ in (\ref{y6}).
 Then the following hold
  \ben
 \item[{\rm (i)}] $\alpha$ and $\beta$  are written as
   \be\label{y11}
 \alpha=4\mu^{-1}(\mu/4\pm\delta^2\mu^{-1}\mp c^2)h,\ \
     \beta=4\mu^{-\frac{3}{2}}\sqrt{\mu/4\pm\delta^2\mu^{-1}\mp c^2}\
     c_0,
  \ee

   \item[{\rm (ii)}] Put $L=L_{\epsilon}^{\pm}$ (for $F=F_{\epsilon}^{\pm}$) the same meaning as that in
   Theorem \ref{th2} (ii). Then
   \be\label{y12}
  L_{\epsilon}^{\pm}=\frac{2\pi}{\sqrt{\mu}}\pm
  \frac{8\pi}{\mu^2\sqrt{\mu}}\ \delta^2.
   \ee

  \item[{\rm (iii)}] Denote by
 $K=K_{\epsilon}^{\pm}$ the scalar flag curvature of $F_{\epsilon}^{\pm}$.
 Then the maximal and minimal values of $K_{\epsilon}^{\pm}$ in three cases
  are as follows:
   \beq
  && Min(K_2^+)=\frac{\big(\sqrt{4\delta^2+\mu^2}-2\delta\big)^3}{\mu\sqrt{4\delta^2+\mu^2}},
   \ \ \
         Max(K_2^+)=\frac{\big(\sqrt{4\delta^2+\mu^2}+2\delta\big)^3}{\mu\sqrt{4\delta^2+\mu^2}},\label{y13}\\
    &&Min(K_0^+)=\frac{\mu^2-8\delta^2}{\mu}, \hspace{2cm}
    Max(K_0^+)=\frac{(\mu^2+4\delta^2)^4}{\mu(\mu^2+8\delta^2)^3},\label{y14}\\
  &&Min(K_0^-)=\frac{\mu^5(\mu^2-16\delta^2)}{(\mu^2-8\delta^2)^3}, \hspace{1.2cm}
    Max(K_0^-)=\frac{(\mu^2-4\delta^2)^4}{\mu(\mu^2-8\delta^2)^3}.\label{y15}
   \eeq

  \een

\end{thm}

By (\ref{y13}) we have the following corollary for square metrics.

\begin{cor}\label{cor4}
 Any  square metric $F=(\alpha+\beta)^2/\alpha$  which is locally
 projectively flat on $S^n$ always has the scalar flag curvature
 with positive lower bound.
\end{cor}

By \cite{SYa} \cite{Y3} and Corollary \ref{cor4}, any  square
metric of scalar flag curvature on $S^n$ ($n\ge 3$) always has the
scalar flag curvature with positive lower bound. Among
$(\alpha,\beta)$-metrics of non-Randers type which are locally
projectively flat on $S^n$, we have not found any other metric
type  having the property for square metrics shown in Corollary
\ref{cor4}. For the metric type $F_0^+$, its scalar flag curvature
has positive lower bound iff. $\mu^2>8\delta^2$ by (\ref{y14}).
For the metric type $F_0^-$ (with the regular condition
$\mu^2>12\delta^2$ by Remark \ref{rem51} below), its scalar flag
curvature has positive lower bound iff. $\mu^2>16\delta^2$ by
(\ref{y15}).

\section{Preliminaries}\label{sec2}

  A geodesic $\gamma=\gamma(t)=(x^i(t))$ of a Finsler metric
  $F=F(x,y)$ are characterized by
   $$\frac{d^2 x^i}{d t^2}+2G^i(x,\frac{d x^i}{d t})=0,\ \ \ \ \
    \Big(G^i:=\frac{1}{4}g^{il}\big \{[F^2]_{x^ky^l}y^k-[F^2]_{x^l}\big \}\Big).$$
In this case, $t$ is the arc-length parameter of $\gamma$. A
Finsler manifold is called complete (positively complete or
negatively complete) if every geodesic $\gamma=\gamma(t)$ is
defined on $(-\infty,+\infty)$ ($(0,+\infty)$ or $(-\infty,0)$)
for the arc-length parameter $t$.

For a Finsler metric $F$, the Riemann curvature $R_y=R^i_{\
k}(y)\frac{\pa}{\pa x^i}\otimes dx^k$ is defined by
 $$
 R^i_{\ k}:=2\frac{\pa G^i}{\pa x^k}-y^j\frac{\pa^2G^i}{\pa x^j\pa
 y^k}+2G^j\frac{\pa^2G^i}{\pa y^j\pa y^k}-\frac{\pa G^i}{\pa y^j}\frac{\pa G^j}{\pa
 y^k}.
 $$
 A Finsler metric $F$ is called of scalar flag
curvature if there is a function $K= K(x,y)$ such that
 $$ R^i_{\ k}=KF^2(\delta^i_k-F^{-2}y^iy_k), \ \ y_k:=(F^2/2)_{y^iy^k}y^i.$$
If $K$ is a constant, then $F$ is called of constant flag
curvature.

Two Finsler metrics $F$ and $\widetilde{F}$ on a same manifold $M$
are called projectively related if they have same geodesics as
point sets, or equivalently their sprays $G^i$ and
$\widetilde{G}^i$ are related by $\widetilde{G}^i=G^i+Py^i$, where
$P=P(x,y)$ is called the projective factor satisfying $P(x,\lambda
y)=\lambda P(x,y)$ for $\lambda>0$. A Finsler metric $F$ is said
to be locally projectively flat if $F$ is projectively related to
a locally Euclidean metric (namely, $G^i=Py^i$ in
 some local coordinate system everywhere on $M$). A locally
 projectively flat Finsler metric ($G^i=Py^i$) is of scalar flag
 curvature $K=K(x,y)$ which is given by
 \be\label{y16}
  K=\frac{P^2-P_{x^k}y^k}{F^2}.
  \ee

An $(\alpha,\beta)$-metric is a Finsler metric defined by a
Riemann metric $\alpha=\sqrt{a_{ij}(x)y^iy^j}$ and a 1-form
$\beta=b_i(x)y^i$ as follows:
 $$F=\alpha \phi(s),\ \ s=\beta/\alpha,$$
 where $\phi(s)>0$ is a $C^{\infty}$ function on
$(-b_o,b_o)$. It is proved in \cite{Shen3} that an
$(\alpha,\beta)$-metric is regular if and only if
 \be\label{y17}
\phi(s)>0,\ \ \ \phi(s)-s\phi'(s)>0,\ \ \
\phi(s)-s\phi'(s)+(b^2-s^2)\phi''(s)>0,\ \ (|s|\le b<b_o).
 \ee
 For an $(\alpha,\beta)$-metric,  the spray coefficients $G^i$ of $F$
 are given by
  \be\label{y18}
  G^i=G^i_{\alpha}+\alpha Q s^i_{\ 0}+\alpha^{-1}\Theta (-2\alpha Q
  s_0+r_{00})y^i+\Psi (-2\alpha Q s_0+r_{00})b^i,
  \ee
where $G^i_{\alpha}$ denote the spray coefficents of $\alpha$ and
 \beqn
  &&Q:=\frac{\phi'}{\phi-s\phi'},\ \
  \Theta:=\frac{Q-sQ'}{2\Delta},\ \
  \Psi:=\frac{Q'}{2\Delta},\ \ \Delta:=1+sQ+(b^2-s^2)Q',\\
  &&r_{ij}:=\frac{1}{2}(b_{i|j}+b_{j|i}),\ \ s_{ij}:=\frac{1}{2}(b_{i|j}-b_{j|i}),\ \
  s_j:=b^is_{ij},\ \ s^i_{\ j}:=a^{im}s_{mj},
 \eeqn
where the covariant derivatives $b_{i|j}$ are taken with respect
to the Levvi-Civita connection of $\alpha$, and we define
$T_0:=T_iy^i$ for a tensor $T=(T_i)$. For an $n$-dimensional
$(\alpha,\beta)$-metric $F=\alpha\phi(\beta/\alpha)$ of
non-Randers type, assume $\beta$ is not parallel with respect to
$\alpha$, and $n\ge 3$ if $F$ is of the metric type $F_0^{\pm}$
(see (\ref{y2})). Then $F$ is locally projectively flat iff.
(\cite{Shen1} \cite{Y1})
  \beq
    &&\big\{1+(k_1+k_3)s^2+k_2s^4\big\}\phi''(s)=(k_1+k_2s^2)\big\{\phi(s)-s\phi'(s)\big\},\label{y19}\\
    && b_{i|j}=\tau
     \big\{(1+k_1b^2)a_{ij}+(k_2b^2+k_3)b_ib_j\big\},\label{y20}\\
     && G^i_{\alpha}=\theta y^i-1/2\cdot\tau (k_1\alpha^2+k_2\beta^2)b^i,\label{y21}
   \eeq
 where  $k_1,k_2,k_3$ are constant with $k_2\ne k_1k_3$, $\theta$ is a 1-form and $\tau=\tau(x)$ is a scalar function. By
 (\ref{y20}) and (\ref{y21}), define a Riemann metric $h=\sqrt{h_{ij}y^iy^j}$ and a
 1-form $\rho_iy^i$ by (\ref{y6}), where $u=u(b^2)>0,v=v(b^2),w=w(b^2)\ne 0$  satisfy
 (\ref{y3})--(\ref{y5}). Then $h$ is of constant sectional curvature
$\mu$ and $\rho$ is closed and conformal with respect to $h$
satisfying (\ref{y7}) (cf. \cite{Y3} \cite{Yu1}).  This fact is
also verified directly in the proof of Theorem \ref{th2} (iii)
(see (\ref{y32}) below). In \cite{Yu1}, it chooses
 \be\label{y22}
   u=e^{2\sigma}, \ \ v=(k_1+k_3+k_2b^2)u, \ \ w=\sqrt{1+(k_1+k_3)b^2+k_2b^4}\ e^{\sigma},
   \ee
   where $\sigma$ is defined by
 $$\sigma:=\frac{1}{2}\int_0^{b^2}\frac{k_2t+k_3}{1+(k_1+k_3)t+k_2t^2}dt.$$
It is easy to verity that (\ref{y22}) satisfies
(\ref{y3})--(\ref{y5}), and in this case,
$||\beta||_{\alpha}=||\rho||_h$, namely, the left hand side of
(\ref{y10}) is just $b^2$.

 \begin{prop}\label{prop21}
  Let $F=\alpha\phi(\beta/\alpha)$ be an $(\alpha,\beta)$-metric
  with $\phi(s)$ being defined by (\ref{y19}) with $k_2\ne
  k_1k_3$ and $\phi(0)=1$. Then $F$ is regular iff. $\phi(s)>0$ for $|s|\le b$ and one of the
  following two cases holds:
  \ben
  \item[{\rm (i)}] If $k_2\le 0$, then $b^2$ is determined by
  $
 1+k_1b^2>0,\ 1+(k_1+k_3)b^2+k_2b^4>0.
  $

\item[{\rm (ii)}] If $k_2> 0$, then $b^2$ is determined by $S\cap
T$, where $S,T$ are two sets defined by
 \beqn
S:=\hspace{-0.6cm}&&\big\{1+k_1b^2>0,\ 1+(k_1+k_3)b^2+k_2b^4>0\},\\
  T:=\hspace{-0.6cm}&&\big\{k_1+k_3\ge 0\big\}\cup\big\{k_1+k_3+2k_2b^2\le
0,1+(k_1+k_3)b^2+k_2b^4>0\} \\
&&\cup\big\{k_1+k_3+2k_2b^2>0,4k_2-(k_1+k_3)^2>0\big\}.
 \eeqn
 \een
\end{prop}

{\it Proof} : By (\ref{y19}) with $\phi(0)=1$, we have
 $
  \phi-s\phi'=exp\big(-\frac{1}{2}\int_0^{s^2}
  \frac{k_1+k_2\theta}{1+(k_1+k_3)\theta+k_2\theta^2}\big)>0,
 $
and
 $$
 \phi-s\phi'+(b^2-s^2)\phi''=(\phi-s\phi')\cdot
 \frac{1+k_1b^2+(k_3+k_2b^2)s^2}{1+(k_1+k_3)s^2+k_2s^4}.
 $$
The above expression is positive for $|s|\le b$ iff.
 $$
 1+k_1b^2+(k_3+k_2b^2)s^2>0,\ \ 1+(k_1+k_3)s^2+k_2s^4>0,
 $$
 from which and the regular condition (\ref{y17}), we can easily
 complete the proof.\qed

\bigskip

By Proposition  \ref{prop21}, the regular condition of the metric
$F_{\epsilon}^{\pm}$ defined by (\ref{y2}) can be determined (for
some special values of $\epsilon$). By putting $k_1=\pm
2,k_2=0,k_3=\mp 3$ in (\ref{y19}), we get the metric
$F_{\epsilon}^{\pm}$. So for $F_{2}^{+}$ and $F_{0}^{+}$, we
obtain the regular condition $b^2<1$, and for $F_{0}^{-}$, it is
$b^2<1/2$.

\section{Proof of Theorem \ref{th1}}

To prove Theorem \ref{th1}, we need to show some lemmas first as
follows.

\begin{lem}\label{lem30} (\cite{IT} \cite{Ta})
 Let $(M,\alpha)$ be an $n$-dimensional complete Riemann manifold with $\alpha=\sqrt{a_{ij}y^iy^j}$ being
 a Riemann metric on the manifold $M$, and $\tau=\tau(x)$ be a
 non-constant scalar function on the manifold $M$. Suppose the covariant
    derivatives $\tau_{i|j}$  with
    respect to $\alpha$ satisfy $\tau_{i|j}=\lambda a_{ij}$ for some scalar function
$\lambda=\lambda(x)$, where $\tau_i:=\tau_{x^i}$. Let
$\gamma=\gamma(t)$ be a geodesic tangent to the gradient vector
field $\nabla \tau$ with $t$ being the arc-length parameter, and
$V$ be a hypersurface (defined by $\tau=constant$) intersecting
with $\gamma$ at a point $P$ with $\nabla\tau(P)\ne 0$. Then
$\tau=\tau(t)$ is a function depending only on $t$ and
$(M,\alpha)$ has a line element in the following form
 $$
 ds^2=(\tau'(t))^2d\widetilde{s}^2+dt^2,
 $$
where $d\widetilde{s}^2$ is a line element of the hypersurface
$V$.
\end{lem}

\begin{lem}\label{lem31}
Let $\alpha=\sqrt{a_{ij}y^iy^j}$ be a a Riemann metric of constant
sectional curvature $\mu$ and $V=V_iy^i$ be a closed 1-form.
Suppose the covariant
    derivatives $V_{i|j}$ of $V$ with
    respect to $\alpha$ satisfy $V_{i|j}=-2ca_{ij}$ for some scalar function
$c=c(x)$. Then we have
 \be\label{y23}
 c_{i|j}+\mu c a_{ij}=0, \ \
 \mu V_i=2c_i, \ \ ||\nabla c||_{\alpha}^2+\mu c^2=constant,
 \ee
 where $\nabla$ is
 the gradient operator with respect to $\alpha$.

\end{lem}

{\it Proof :} Part of (\ref{y23}) has been proved (\cite{KN}). We
also show here. By $V_{i|j}=-2ca_{ij}$ and the Ricci identity we
have
 $$
2(c_ka_{ij}-c_ia_{jk})-V_mR^{\ m}_{j\ ki}=0,
 $$
where $R$ is the Riemann curvature tensor of $\alpha$. Since
$\alpha$ is of constant sectional curvature $\mu$, the above
equations show
 $$
 2(c_ka_{ij}-c_ia_{jk})=\mu(V_ka_{ij}-V_ia_{jk}),
 $$
which are equivalent to
 $$(2c_k-\mu V_k)a_{ij}=(2c_i-\mu V_i)a_{jk}.$$
Contracting the above by $a^{jk}$ we obtain $2c_i=\mu V_i$. Thus
we get the former two formula in (\ref{y23}).

 Now by the first formula in (\ref{y23})
we have
 $$(c^ic_i+\mu c^2)_{|k}=2(c^ic_{i|k}+\mu cc_k)=2(-\mu cc_k+\mu cc_k)=0,$$
 where $c^i:=a^{ij}c_j$. So $||\nabla c||_{\alpha}^2+\mu c^2$ is a
 constant.     \qed

\begin{lem}\label{lem32}
 Let $(M,F)$ and $(M,\widetilde{F})$ be two Finsler manifolds with
 the metric relation $\widetilde{F}\le F$. If
 $(M,\widetilde{F})$ is complete (resp. positively, or negatively
 complete), then $(M,F)$ is also complete (resp. positively, or negatively
 complete).
\end{lem}

{\it Proof }: Let $(M,\widetilde{F})$ be positively complete. It
is easy to see that any forwarded Cauchy sequence induced by the
metric $d_F$ is also a forwarded Cauchy sequence induced by the
metric $d_{\widetilde{F}}$. Then we complete the proof by
Hopf-Rinow Theorem for Finsler manifolds. \qed

\begin{lem}\label{lem33}
 Let $F=\alpha\phi(\beta/\alpha)$ satisfy the condition  (C5)
 in Theorem \ref{th1}, where
 $\phi(s)$ is determined by (\ref{y19}) with $k_2\ne k_1k_3$.
 Define $h$ and $\rho$ by (\ref{y6}), where $u,v,w$ is given by
 (\ref{y22}). Then $F$ is positively complete iff. $\alpha$ is complete iff. $h$ is
 complete.
\end{lem}

{\it Proof} : By the condition (C5) $b(x_0)=Sup_{x\in M} b(x)$ at
some point $x_0\in M$, and the continuity of the positive function
$\phi(s)$ on $|s|\le b(x_0)$, we get two constants $m_1>0$ and
$M_1>0$ such that $m_1\alpha\le F\le M_1\alpha$. So $F$ is
positively complete iff. $\alpha$ is complete by Lemma
\ref{lem32}.

By the choice of $u,v,w$ given by (\ref{y22}), we have
 $$
 1+\frac{v}{u}\frac{\beta^2}{\alpha^2}=1+\frac{v}{u}s^2=1+(k_1+k_3+k_2b^2)s^2.
 $$
Therefore we get
 \beqn
 &&k_1+k_3+k_2b^2\ge 0:\ \ \  1\le 1+\frac{v}{u}\frac{\beta^2}{\alpha^2}\le
 1+(k_1+k_3+k_2b^2)b^2,\\
&&k_1+k_3+k_2b^2< 0: \ \ \ 1+(k_1+k_3+k_2b^2)b^2\le
1+\frac{v}{u}\frac{\beta^2}{\alpha^2}\le 1.
 \eeqn
Now put
 $$A_1(b^2):=Min\big\{1,1+(k_1+k_3+k_2b^2)b^2\big\},\ \ \
 A_2(b^2):=Max\big\{1,1+(k_1+k_3+k_2b^2)b^2\big\}.$$
Then it follows from the regular condition of $F$ shown in
Proposition \ref{prop21} that both $A_1(t)$ and $A_2(t)$ are
positive and continuous functions on the closed interval $|t|\le
b(x_0)$. Let
 $$
 m_2:=Min\big(\sqrt{u(t)A_1(t)},\ |t|\le b(x_0)\big),\ \ M_2:=Max\big(\sqrt{u(t)A_2(t)},\ |t|\le b(x_0)\big)
 $$
 Then by continuity, we have $M_2\ge m_2>0$. Therefore, we obtain
  $$m_2\alpha\le h=\sqrt{u\big(1+\frac{v}{u}\frac{\beta^2}{\alpha^2}}\big)\cdot\alpha\le
  M_2\alpha,
  $$
which implies that $\alpha$ is complete iff. $h$ is complete by
Lemma \ref{lem32}.  \qed

\begin{lem}\label{lem34}
 Let $F=\alpha\phi(\beta/\alpha)$ satisfy
 (\ref{y19})--(\ref{y21}). Then $F$ is projectively related to the
 Riemann metric $h$ defined by (\ref{y6}), where $u,v$ are
 determined by (\ref{y3})--(\ref{y5}).
\end{lem}

{\it Proof :} In some local coordinate, the spray $G^i_{\alpha}$
are given by (\ref{y21}). Plugging (\ref{y20}), (\ref{y21}) and
$$s_0=0,\ \ s^i_{\ 0}=0,\ \ \Psi=\frac{1}{2}\frac{\phi''}{\phi-s\phi'+(b^2-s^2)\phi''}$$
into (\ref{y18}), and then using (\ref{y19}), we easily get the
spray $G^i$ of $F$ in the form $G^i=P_1y^i$ for some positively
homogeneous function $P_1$ of degree one on $TM$. Now define a
Riemann metric $h$ shown in (\ref{y6}). In the same local
coordinate as shown in (\ref{y21}) for $G^i_{\alpha}$, using
(\ref{y20}), (\ref{y21}) and (\ref{y3})--(\ref{y5}), one direct
computation shows that the spray $G^i_h$ of $h$ satisfy
$G^i_h=P_2y^i$ for some positively homogeneous function $P_2$ of
degree one on $TM$. Therefore, the sprays $G^i$ and $G^i_h$
satisfy the relation $G^i=G^i_h+Py^i$ for $P:=P_1-P_2$ under any
local coordinate. Thus $F$ and $h$ are projectively related. \qed

\bigskip

Now we start the proof of Theorem \ref{th1}. Let $(M,F)$ be an
$n(\ge 2)$-dimensional Finsler manifold, where
$F=\alpha\phi(\beta/\alpha)$ is an $(\alpha,\beta)$-metric on $M$.
If $\beta$ is parallel with respect to $\alpha$, then by
(\ref{y18}) we get $G^i=G^i_{\alpha}$. Thus $\alpha$ is of
constant sectional curvature  since $F$ is locally projectively
flat. If $\beta=0$, then $F=\alpha$ is Riemannian. If $\beta\ne
0$, then it is easy to show that $\alpha$ is flat, and so $F$ is
flat-parallel.

Suppose $\beta$ is not parallel with respect to $\alpha$. Then by
the conditions (C1), (C2) and (C4), the $(\alpha,\beta)$-metric
$F=\alpha\phi(\beta/\alpha)$ satisfies (\ref{y19})--(\ref{y21})
with $k_2\ne k_1k_3$ (see \cite{Shen1} \cite{Y1}). Now make a
transformation (\ref{y6}), namely,
 $$
h:=\sqrt{u\alpha^2+v\beta^2},\ \ \ \rho:=w\beta,
 $$
where $u=u(b^2)>0,v=v(b^2),w=w(b^2)\ne 0$ are suitable functions
satisfying the ODEs (\ref{y3})--(\ref{y5}). As shown in Section
\ref{sec2}, $h=\sqrt{h_{ij}y^iy^j}$ is a Riemann metric of
constant sectional curvature $\mu$ and the 1-form $\rho=p_iy^i$ is
closed and conformal with respect to $h$ (cf. \cite{Y3}
\cite{Yu1}). So the covariant    derivatives $p_{i|j}$ of $\rho$
with
    respect to $h$ satisfy (\ref{y7}), that is,
 \be\label{y24}
 p_{i|j}=-2ch_{ij}
 \ee
 for some scalar function $c=c(x)$. Then by Lemma \ref{lem31}, the
 scalar function $c$ in (\ref{y24}) satisfies
  \be\label{y25}
c_{i|j}+\mu c h_{ij}=0.
  \ee
  Now we fix a special choice of $u,v,w$ given by
  (\ref{y22}), namely,
 \beqn
   &&u=e^{2\sigma}, \ \ v=(k_1+k_3+k_2b^2)u, \ \ w=\sqrt{1+(k_1+k_3)b^2+k_2b^4}\
   e^{\sigma},\\
 &&\sigma:=\frac{1}{2}\int_0^{b^2}\frac{k_2t+k_3}{1+(k_1+k_3)t+k_2t^2}dt.
 \eeqn
Note that in this case, $b^2=||\beta||_{\alpha}^2=||\rho||_h^2$,
and  $h$ is complete by Lemma \ref{lem33}.

\

\noindent {\bf Case I:} Assume $c=constant$ on the manifold $M$.
We will prove $c=0$ in this case.

Suppose $c\ne 0$. Then by (\ref{y25}) we get $\mu=0$, that is, $h$
is locally Euclidean. Since $(M,h)$ is complete by Lemma
\ref{lem33}, its universal covering metric space is
$(R^n,\widetilde{h})$ which is locally isometric with $(M,h)$,
where $\widetilde{h}=|y|$. The lift $\widetilde{\rho}$ of $\rho$
in $(R^n,\widetilde{h})$ satisfies
$\widetilde{p}_{i|j}=-2c\widetilde{h}_{ij}$ (see (\ref{y24})),
from which it is easy to get
$\widetilde{p}_i=\widetilde{p}^i=-2cx^i+\xi^i$ under the global
coordinate $(x^i)$ of $\widetilde{h}$, where $\xi$ is a constant
vector. But $b^2=||\beta||_{\alpha}^2=||\rho||_h^2\le b^2(x_0)$ is
bounded, so
 $$||\rho||_h^2=||\widetilde{\rho}||_{\widetilde{h}}^2=4c^2|x|^2-4c\langle
 \xi,x\rangle+|\xi|^2$$
 is also bounded on $R^n$. However, it is unbounded since $c\ne 0$, which is a
 contradiction. Thus we get $c=0$ if $c=constant$ on the manifold $M$.

 Now by $c=0$, it follows from (\ref{y24}) that $\rho$ is parallel
 with respect to $h$.
 If  $\rho=0$, it is clear that $F=\alpha$. If $\rho\ne 0$, then $h$ is locally
 flat. Now by (\ref{y6}) we get
 $$
 \alpha^2=\frac{1}{u}h^2-\frac{v}{uw^2}\rho^2,\ \
 \beta=\frac{1}{w}\rho.
 $$
The above shows that $\alpha$ is flat and $\beta$ is parallel with
respect to $\alpha$, since locally $\alpha$ and $\beta$ are
independent of $x\in M$, which follows from the fact that
$u=u(b^2),v=v(b^2),w=w(b^2)$ are constant and locally $h=|y|$ and
$\rho=\langle \xi,y\rangle$ for some constant vector $\xi$.

\

\noindent {\bf Case II:} Assume $c$ is a non-constant. In this
case, we will show that $\mu>0$ and the manifold $M$ is
homeomorphic to the $n$-sphere $S^n$.

If $\mu=0$, then by (\ref{y23}) in Lemma \ref{lem31}, we have
$2c_i=\mu p_i$. So $c_i=0$, which shows that $c$ is a constant.
This is a contradiction.

Now assume $\mu\ne 0$. It is clear that (\ref{y25}) implies that
the trajectories of the gradient vector field $\nabla c$ are
geodesics of $h$ as point sets. Let $x(t)$ be a geodesic of
$(M,h)$ tangent to $\nabla c$, where $t$ is the arc-length
parameter. Put $c(t)=c(x(t))$. Then by (\ref{y25}) we have
 \be\label{y26}
 \frac{d^2c(t)}{dt^2}+\mu c(t)=0.
 \ee

If $\mu<0$, then the solutions of (\ref{y26}) are given by
 \be\label{y27}
 c(t)=ke^{\sqrt{-\mu}\ t}+le^{-\sqrt{-\mu}\ t},
 \ee
where $k,l$ are constant with $k^2+l^2\ne 0$. Since $(M,h)$ is
complete by Lemma \ref{lem33}, $c(t)$ in (\ref{y27}) is defined on
$ (-\infty,+\infty)$. By (\ref{y23}) in Lemma \ref{lem31}, we have
$2c_i=\mu p_i$. By Lemma \ref{lem30}, we get $||\nabla
c||_h^2=(c'(t))^2$. So $b^2=||\rho||_h^2=4\mu^{-2}(c'(t))^2$ is
unbounded for $t\in (-\infty,+\infty)$ by (\ref{y27}). But by the
condition (C5) in Theorem \ref{th1}, $b^2$ is bounded on $M$. We
get a contradiction.

The above discussion shows that there must have $\mu>0$. Since $c$
is a non-constant function satisfying (\ref{y25}) and $(M,h)$ is
complete, it follows from a result in \cite{Ob} \cite{Ta} that the
manifold $M$ is homeomorphic to the $n$-sphere $S^n$.

\

The discussions in Case I and Case II above have actually proved
the following lemma.
\begin{lem}\label{lem36}
 Let $h=\sqrt{h_{ij}y^iy^j}$ be a complete
Riemannian metric of constant sectional curvature $\mu$ and the
1-form $\rho=p_iy^i$ be closed and conformal with respect to $h$
satisfying (\ref{y24}) for some scalar function $c=c(x)$. Suppose
$||\rho||_h$ is bounded. Then $\mu=0$  and $\rho$ is parallel with
respect to $h$, or $\mu>0$ and $M$ is homeomorphic to the
$n$-sphere $S^n$.
\end{lem}

Now we prove the final part in Theorem \ref{th1}, namely, if $F$
satisfies the conditions (C3)--(C5), and
 (C6) $F$ is of constant flag curvature, then $F=\alpha$ or $F$ is
 flat-parallel. According to the results in \cite{LS} \cite{Y2},
 if an $(\alpha,\beta)$-metric $F=\alpha\phi(\beta/\alpha)$ is
 locally projectively flat with constant flag curvature, then $F$
 is one of the three classes: (a1) $F$ is flat-parallel; (a2) $F=\alpha+\beta$ is of
 Randers type; (a3) $F=(\alpha+\beta)^2/\alpha$ is of square type.

For a Randers metric $F=\alpha+\beta$ (the regular condition $b^2<
1$) which is locally projectively flat with constant flag
curvature, it follows from \cite{Shen0} that $\alpha$ is of
non-positive constant sectional curvature $\mu\le 0$ and the
covariant derivative of $\beta$ with respect to $\alpha$ satisfies
  \be\label{y28}
 b_{i;j}=2k (a_{ij}-b_ib_j),
 \ee
 for some constant $k$. Define a Riemann metric
 $h=\sqrt{h_{ij}y^iy^j}$ and 1-form $\rho=p_iy^i$ by
  $$
 h:=\alpha,\ \ \ \rho:=\frac{1}{\sqrt{1-b^2}}\beta.
  $$
Under the above transformation, it can be directly verified that
$||\rho||_h^2=b^2/(1-b^2)$, and
  (\ref{y28}) is equivalent to
  $$
 p_{i|j}=-2ch_{ij}, \ \ \ \big(c:=-\frac{k}{\sqrt{1-b^2}}\big),
   $$
where the covariant derivative is taken
  with respect to $h=\alpha$.
  Thus (\ref{y24})  holds for some scalar function $c=c(x)$. Further, since $F=\alpha+\beta\le
  2\alpha$, it follows from Lemma \ref{lem32} that $h=\alpha$ is
  complete. Now it follows from  the condition (C5)  and $||\rho||_h^2=b^2/(1-b^2)$ that
 $||\rho||_h$ is bounded. Plus  the fact $\mu\le0$,
  it directly follows from Lemma \ref{lem36} that $F=\alpha$ or $F$ is
  flat-parallel.

For a square metric $F=(\alpha+\beta)^2/\alpha$ (the regular
condition $b^2< 1$)  which is locally projectively flat with
constant flag curvature, consider a special transformation of
(\ref{y6}) to define a new pair $(h,\rho)$ by putting
$u=(1-b)^2,v=0$, and $w=\sqrt{1-b^2}$, namely,
 $$
 h:=(1-b^2)\alpha, \ \ \ \rho:=\sqrt{1-b^2}\beta.
 $$
it is already shown that $h$ is of constant sectional curvature
$\mu$ and $\rho$ is a closed and conformal 1-form with respect to
$h$ satisfying (\ref{y24}) for some scalar function $c=c(x)$. In
\cite{SY}, it actually proves $\mu\le 0$ (also see \cite{Y3}).
Then it follows from the condition (C5)  and
$||\rho||_h^2=b^2/(1-b^2)$ that
 $||\rho||_h$ is bounded. Plus  the fact $\mu\le0$,
  it directly follows from Lemma \ref{lem36} that $F=\alpha$ or $F$ is
  flat-parallel. \qed

\begin{rem}\label{rem31}
 Let $\phi(s)$ be a fixed function with $|s|<b_o$ (possibly $b_o=+\infty$).
 Let $\hat{b}\le b_o$ be a constant such that
 for arbitrary $\alpha$ and $\beta$ with
 $k<||\beta||_{\alpha}<\hat{b}$ (where $k$ is a constant with $0\le
 k<\hat{b}$),
 $F=\alpha\phi(\beta/\alpha)$ is always regular, and on the
 other hand, if a pair of $\alpha$ and $\beta$ satisfies $||\beta||_{\alpha}\ge
 \hat{b}$, then $F=\alpha\phi(\beta/\alpha)$ is non-regular. Note
 that $\hat{b}$ is dependent on the function $\phi$ and
 independent of $\alpha,\beta$. For example, if $\phi(s)=1+s$ or
 $(1+s)^2$ or $1+s^2$, we have $\hat{b}=1$; if $\phi(s)=1-s^2$, we have
 $\hat{b}=1/\sqrt{2}$.

 Now
 in Theorem \ref{th1}, the condition (C5) can be replaced by a
 weaker condition $b=||\beta||_{\alpha}\le \bar{b}<\hat{b}$ for a
 constant $\bar{b}$.
\end{rem}

\section{Proof of Theorem \ref{th2}}

Let $\phi(s)$ satisfy the ODE (\ref{y1}), and
$F=\alpha\phi(\beta/\alpha)$ is an $(\alpha,\beta)$-metric on
$S^n$. Assume the conditions  (C2) and (C4) in Theorem \ref{th1}
hold. Then $\beta$ is not parallel with respect to $\alpha$.
Otherwise, $\alpha$ is flat on $S^n$ since $\beta\ne 0$ ($F$ is
non-Riemannian), which is impossible. Now as shown in the proof of
Theorem \ref{th1}, define a pair $(h,\rho)$ (see (\ref{y6})) for
every suitable $u,v,w$ satisfying (\ref{y3})--(\ref{y5}). Then
$h=\sqrt{h_{ij}y^iy^j}$ is a Riemann metric of constant sectional
curvature $\mu>0$ and $\rho=p_iy^i$ is a closed and conformal
1-form with respect to $h$ satisfying (\ref{y24}) and then
(\ref{y25}) for a scalar function $c=c(x)$.

By Lemma \ref{lem31}, we have $\rho=2\mu^{-1}c_0$ and
$\delta:=\sqrt{||\nabla c||_h^2+\mu c^2}$ is a constant.  Put
$c(t)=c(x(t))$, where $x(t)$ is a geodesic of $(M,h)$ tangent to
$\nabla c$ with $t$ being the arc-length parameter. Solving the
ODE (\ref{y26}) for $\mu>0$ we obtain $c=\mu^{-\frac{1}{2}}\delta
\cos(\sqrt{\mu}\ t)$ for a suitably chosen again  arc-length
parameter $t$. By Lemma \ref{lem30}, we have $||\nabla
c||_h=|c'(t)|$. Thus the gradient field $\nabla c$ has just two
vanishing points $P,Q\in S^n$. This completes the proof of item
(i) of Theorem \ref{th2}.

To prove item (ii) of Theorem \ref{th2}, we choose $u,v,w$ given
by (\ref{y22}). By Lemma \ref{lem34}, the $(\alpha,\beta)$-metric
$F$ is projectively related to the Riemann metric $h$. So $F$ and
$h$ have same geodesics as point sets. In the following we
consider the $F$-length of two families of geodesics of $F$.

The two vanishing points $P,Q$ of the gradient field $\nabla c$
are just a pair of antipodal points on $S^n$. Let $x(t)$ be an
arbitrary closed geodesic of $h$ connecting $P$ and $Q$ (great
circle of $h$) with $0\le t\le 2\pi/\sqrt{\mu}$, and
$c=\mu^{-\frac{1}{2}}\delta \cos(\sqrt{\mu}\ t)$ along $x(t)$.
From (\ref{y6}) and $\rho=2\mu^{-1}c_0$, we have
 \be\label{y29}
 \alpha^2=u^{-1}\big[h_{\mu}^2-4(\mu w)^{-2}vc_0^2\big],\ \
     \beta=2(\mu w)^{-1} c_0.
 \ee
Then along the closed geodesic $x(t)$ of $h$, from
$h(x(t),x'(t))=1$, $c_0(x(t),x'(t))=c'(t)$,
$b^2=||\rho||_h^2=4\mu^{-2}(c'(t))^2$ and the choice of $u,v,w$
defined by (\ref{y22}), we obtain from (\ref{y29})
 $$
 \alpha(x(t),x'(t))=\frac{1}{w}\sqrt{\frac{w^2-b^2v}{u}}=\frac{1}{w}=
 \frac{e^{-\frac{1}{2}\int_0^{b^2}\frac{k_2\theta+k_3}{1+(k_1+k_3)\theta+k_2\theta^2}d\theta}}{\sqrt{1+(k_1+k_3)b^2+k_2b^4}},
 $$
 $$
 \frac{\beta(x(t),x'(t))}{\alpha(x(t),x'(t))}=\frac{2\mu^{-1}w^{-1}\sqrt{u}\ c'(t)}{\sqrt{1-4\mu^{-2}w^{-2}v(c'(t))^2}}
 =\frac{2\mu^{-1}c'(t)}{\sqrt{u^{-1}(w^2-b^2v)}}=2\mu^{-1}c'(t).
 $$
Therefore, the $F$-length $L_1$ of the geodesic $x(t)$ is given by
 \be\label{y30}
     L_1=\int_0^{\frac{2\pi}{\sqrt{\mu}}}\alpha\big(x(t),x'(t)\big)\phi\Big(\frac{\beta(x(t),x'(t))}{\alpha(x(t),x'(t))}\Big)
     =\int_0^{\frac{2\pi}{\sqrt{\mu}}} \frac{\phi\big(-2\mu^{-1}\delta \sin(\sqrt{\mu}\ t)\big)}
     {\sqrt{1+(k_1+k_3)B+k_2B^2}\ e^{\sigma}}dt,
     \ee
  where
   $$B:=4\mu^{-2}\delta^2\sin^2(\sqrt{\mu}\ t)),\ \
   \sigma:=\frac{1}{2}\int_0^B\frac{k_2\theta+k_3}{1+(k_1+k_3)\theta+k_2\theta^2}d\theta,
   $$
By the ODE (\ref{y19}) on $\phi(s)$, we easily get
 $$
 \phi(s)=1+\phi'(0)s+\frac{1}{2}k_1s^2+\big(\frac{1}{12}k_2-\frac{1}{8}k_1^2-\frac{1}{12}k_1k_3\big)s^4+o(s^6).
 $$
Using the above expansion on $\phi(s)$, expanding (\ref{y30}) as
to $\delta$ at $\delta=0$ we obtain (\ref{y8}).

Now we compute the $F$-length  of the family of closed geodesics
lying on the hypersurface $c=0$. It is clear that the hypersurface
$c=0$ is totally geodesic in $(S^n,h)$. Let $z(t)$ is a closed
geodesic of $(S^n,h)$ lying on the hypersurface $c=0$, where $t$
is the arc-length parameter with respect to $h$. Then $z(t)$ is
also a geodesic of $(S^n,F)$ as a point set. It is easily seen
that $c_0(z(t),z'(t))=0$, and $b^2(z(t))=4\mu^{-2}\delta^2$. Thus
it follows from (\ref{y29}) and then from (\ref{y22}) that the
$F$-length $L_2$ of $z(t)$ is given by
 \beq\label{y31}
 L_2&=&\int_0^{\frac{2\pi}{\sqrt{\mu}}}\frac{1}{\sqrt{u}}\ dt=
 \int_0^{\frac{2\pi}{\sqrt{\mu}}}\exp\Big(-\frac{1}{2}\int_0^{\frac{4\delta^2}{\mu^2}}
 \frac{k_2\theta+k_3}{1+(k_1+k_3)\theta+k_2\theta^2}d\theta\Big)\ dt,\nonumber\\
 &=&\frac{2\pi}{\sqrt{\mu}}\exp\Big(-\frac{1}{2}\int_0^{\frac{4\delta^2}{\mu^2}}
 \frac{k_2\theta+k_3}{1+(k_1+k_3)\theta+k_2\theta^2}d\theta\Big).
  \eeq
Now expanding (\ref{y31}) as to $\delta$ at $\delta=0$ we also
obtain (\ref{y8}).

\begin{rem}\label{rem41}
 The two integrals given by (\ref{y30}) and (\ref{y31}) should be
 equal. But we have not found a direct way to prove it.
\end{rem}

To show the estimation in (iii) for the maximal and minimal values
of the scalar flag curvature $K$, we first compute the expression
of $K$.  Since $F$ defined by (\ref{y19}) with $k_2\ne k_1k_3$
satisfies the conditions (C2) and (C4) in Theorem \ref{th1}, we
have (\ref{y20}) and (\ref{y21}). Define a Riemann metric
$h=\sqrt{h_{ij}y^iy^j}$ and 1-form $\rho=p_iy^i$ by (\ref{y6}),
where $u=u(b^2)>0,v=v(b^2),w=w(b^2)\ne 0$ are arbitrary suitable
functions satisfying the ODEs (\ref{y3})--(\ref{y5}). Now using
(\ref{y3})--(\ref{y5}), it follows from (\ref{y20}) and
(\ref{y21}) that
 \be\label{y32}
 G^i_h=\Big(\frac{(v-k_1u)\tau}{u}\beta+\theta\Big)y^i,\ \ \
 p_{i|j}=\frac{w\tau}{u}\ h_{ij}=-2ch_{ij},
 \ee
where the covariant derivatives are taken with respect to $h$. So
by (\ref{y32}), $h$ is locally projectively flat (equivalently,
$h$ is of constant sectional curvature $\mu$) and $\rho$ is closed
and conformal with respect to $h$ (cf. \cite{Y3} \cite{Yu1}). In
some local coordinate, $h$ and its spray are given by
 \be\label{y33}
 h=\frac{\sqrt{(1+\mu |x|^2)|y|^2-\mu\langle
 x,y\rangle^2}}{1+\mu|x|^2},\ \ \ G^i_h=-\frac{\mu \langle x,y\rangle}{1+\mu
 |x|^2}\ y^i.
 \ee
Then the function $c$ in (\ref{y32}) is locally given by
 \be\label{y34}
 c=\frac{-k+\mu\langle \xi,x\rangle}{2\sqrt{1+\mu|x|^2}},
 \ee
where $k$ is a constant and $\xi$ is a vector. By
(\ref{y32})--(\ref{y34}) we obtain
 \be\label{y35}
 \theta=\frac{(k_1u-v)\tau}{u}\beta-\frac{\mu \langle x,y\rangle}{1+\mu |x|^2},\ \ \
 \ \tau=\frac{k-\mu\langle
 \xi,x\rangle}{\sqrt{1+\mu|x|^2}}\frac{u}{w}.
 \ee
Put $\theta_0=\theta_iy^i=\theta_{x^i}y^i$ and
$\tau_0=\tau_iy^i=\tau_{x^i}y^i$. Then by (\ref{y35}) we get
 \beq
 \theta_0&=&\frac{k_1u-v}{u}\big(\tau^2\alpha^2+\tau_0\beta+2\tau\theta\beta\big)+\Big\{k_1k_3-2k_2+(k_3+2k_1)\frac{v}{u}
 -\frac{2v^2}{u^2}\Big\}\tau^2\beta^2\nonumber\\
 &&-\frac{\mu[(1+\mu|x|^2)|y^2|-2\mu\langle
 x,y\rangle^2]}{(1+\mu|x|^2)^2},\label{y36}\\
 \tau_0&=&\frac{\tau(k_3u-v)(k-\mu\langle
 \xi,x\rangle)}{w\sqrt{1+\mu|x|^2}}\beta-\frac{\mu
 u}{\sqrt{1+\mu|x|^2}}\Big\{\frac{(k-\mu\langle
 \xi,x\rangle)\langle x,y\rangle}{w(1+\mu|x|^2)}+\langle
 \xi,y\rangle\Big\},\label{y37}
 \eeq
where we have used (\ref{y20}), (\ref{y21}), the ODEs
(\ref{y3})--(\ref{y5}) and
 $$u_{x^i}y^i=2u'(b^2)r_0,\ \ v_{x^i}y^i=2v'(b^2)r_0,\ \ \beta_{x^i}y^i=r_{00}+2b_mG_{\alpha}^m.$$
Plugging (\ref{y19})--(\ref{y21}) into (\ref{y18}) we obtain
$G^i=Py^i$, where $P$ is given by (\cite{LS})
 $$
 P=\theta+\frac{1}{2}\tau\Big\{\frac{1+(k_1+k_3)s^2+k_2s^4}{\phi}\phi'-(k_1+k_2s^2)s\Big\}\alpha,
 $$
and then by (\ref{y16}), we obtain the scalar flag curvature $K$
given by
 \beq\label{y38}
 K=\hspace{-0.5cm}&&\frac{\theta^2-\theta_0}{\alpha^2\phi^2}+\frac{1}{2\phi^2}\Big\{(k_1+k_2s^2)s
 -\big[1+(k_1+k_3)s^2+k_2s^4\big]\frac{\phi'}{\phi}\Big\}\frac{\tau_0}{\alpha}\nonumber\\
 &&-\frac{\big[1+(k_1+k_3)s^2+k_2s^4\big](k_1+2k_3+3k_2s^2)s\tau^2\phi'}{2\phi^3}
 +\big[1+(k_1+k_3)s^2+k_2s^4\big]^2\times\nonumber\\
 &&\frac{3\tau^2(\phi')^2}{4\phi^4}+\frac{\big[4k_2-k_1^2+2k_2(k_1+2k_3)s^2+3k_2^2s^4\big]s^2\tau^2}{4\phi^2}.
 \eeq
 By Lemma \ref{lem31} we have $2c_0=\mu\rho$
(see (\ref{y23})), and so $2c_0=\mu w\beta$. Using this fact,
differentiating both sides of (\ref{y34}) with respect to $x^i$
gives
 \be\label{y39}
 (k-\mu\langle \xi,x\rangle)\langle
 x,y\rangle+(1+\mu|x|^2) \big[\langle\xi,y\rangle-w\sqrt{1+\mu|x|^2}\ \beta\big]=0.
 \ee
By $h=\sqrt{u\alpha^2+v\beta^2}$ and (\ref{y33}) we get
 \be\label{y40}
 \frac{(1+\mu |x|^2)|y|^2-\mu\langle
 x,y\rangle^2}{(1+\mu|x|^2)^2}=u\alpha^2+v\beta^2.
 \ee
Now plugging the expressions of $\theta,\tau$ and
$\theta_0,\tau_0$ given by (\ref{y35})--(\ref{y37}) into the
scalar flag curvature $K$ shown in (\ref{y38}), we can obtain the
local expression of $K$ under a local coordinate system shown in
(\ref{y33}) for $h$. Then using (\ref{y34}), (\ref{y39}) and
(\ref{y40}), the scalar flag curvature $K$ can be further written
in the form $R(\beta/\alpha,c^2)$, where $R(s,t)$ is a function of
two variables $(s,t)$ defined by (\ref{y9}).

Next we show that $R(s,t)$ is defined on a bounded and closed
subset  in the Euclidean plane.  Since $\rho=2\mu^{-1}c_0$ by
Lemma \ref{lem31}, we have
$$
||\rho||_h^2=4\mu^{-2}||\nabla c||_h^2=4\mu^{-2}(\delta^2-\mu
c^2).
$$
Since $h$ and $\rho$ are defined by $\alpha$ and $\beta$ in
(\ref{y6}), we easily get
 $$
||\rho||_h^2=\frac{w^2(b^2)b^2}{u(b^2)+v(b^2)b^2}.
 $$
Thus we obtain (\ref{y10}), namely,
 \be\label{y41}
 \frac{w^2(B)B}{u(B)+v(B)B}=\frac{4(\delta^2-\mu
       t)}{\mu^2},
 \ee
where $B:=b^2$ and $t:=c^2$. Note that the function on the left
hand side of (\ref{y41}) is strictly increasing on the variable
$B$, which follows from
 $$
 \frac{d}{dB}\Big(\frac{w^2(B)B}{u(B)+v(B)B}\Big)=\frac{w^2}{u\big[1+(k_1+k_3)B+k_2B^2\big]}>0
 $$
because of the ODEs (\ref{y3})--(\ref{y5}) and
$1+(k_1+k_3)b^2+k_2b^4>0$ by the regular condition shown in
Proposition \ref{prop21}. So in the equation (\ref{y41}), $B$ is
uniquely determined for every $t$, and thus $B=B(t)$ is a function
of $t$. Now since $0\le c^2\le \delta^2/\mu$ by Theorem \ref{th2}
(i) and $|\beta/\alpha|\le b$, we see $R(s,t)$ is defined on the
bounded and closed subset $D=\{(s,t)|0\leq t\leq \delta^2/\mu,\
s^2\leq B\}$ in the Euclidean plane, where $B=B(t)$ is determined
by the equation (\ref{y41}).  \qed

\begin{rem}\label{rem42}
 In Theorem \ref{th2} (iii), we can rewrite the function
 $R=R(s,t)$ of two variables as  a different function
  $\widetilde{R}=\widetilde{R}(s,t)$ of two variables. Put
  \beq\label{y42}
 \widetilde{R}=&&\hspace{-0.5cm}\frac{u}{\phi^2
 w^2}\Big\{3u\big(\frac{f_1^2\phi'^2}{\phi^2}+\frac{4f_1-f_2^2}{s^2}\big)
 +2(uf_3-v)\big(f_2-\frac{f_1s\phi'}{\phi}\big)\Big\}A+
 \frac{\mu u (f_1s\phi'-f_2\phi)}{2\phi^3},\\
&&\hspace{2.5cm}
\Big(A:=\frac{1}{\mu}\big(\delta^2-\frac{\mu^2}{4}\frac{w^2
t^2}{u+v t^2}\big)\ \ by \  (\ref{y41})\Big),\nonumber
  \eeq
  where $u,v,w$ are functions of $t^2$ and $(s,t)\in \widetilde{D}$ with
  $\widetilde{D}=\{(s,t)|0\le t\le t_o,|s|\le t\}$, in which $t_o$ is the
  unique positive constant satisfying
   $$
 \frac{w^2(t_o^2)t_o^2}{u(t_o^2)+v(t_o^2)t_o^2}=\frac{4\delta^2}{\mu^2},
   $$
\end{rem}

\section{Proof of Theorem \ref{th3}}

 Let $(S^n,F)$ be an  $n(\ge 2)$-dimensional non-Riemannian Finsler
 manifold with
$F=F_{\epsilon}^{\pm}$ defined by (\ref{y2}), and the conditions
(C2) and (C4) in Theorem \ref{th1} be satisfied. As shown in the
proof of Theorem \ref{th2}, $\beta$ is not parallel with respect
to $\alpha$. Since $\phi(s)=1+\epsilon s\pm s^2$ for the metric
$F=F_{\epsilon}^{\pm}$ defined in (\ref{y2}), we may put $k_1=\pm
2,k_2=0,k_3=\mp 3$ in (\ref{y3})--(\ref{y5}). So in this case, we
may choose
 \be\label{y43}
u=(1\mp B)^2,\ v=0, \ w=\sqrt{1\mp  B},  \ \ \ (B:=b^2).
 \ee
Then using (\ref{y43}),  we get a pair $(h,\rho)$ by (\ref{y6}),
where $h$ is a Riemann metric of constant sectional curvature
$\mu>0$ and $\rho$ is a closed and conformal 1-form.
Correspondingly we get a scalar function $c=c(x)$ satisfying
(\ref{y7}) and a positive constant $\delta$ defined by  $c$. By
the choice of $u,v,w$ shown in (\ref{y43}), it follows from
(\ref{y41}) (therein $t=c^2$) that
 \be\label{y44}
 B=\frac{4(\delta^2-\mu c^2)}{\mu^2\pm
4(\delta^2-\mu c^2) }.
 \ee

Now for the proof of item (i) in this theorem, using (\ref{y43})
and (\ref{y44}), it is clear from (\ref{y6}) that $\alpha$ and
$\beta$ can be expressed by $h,\mu,c,\delta$ shown in (\ref{y11}),
namely,
 \be\label{y45}
 \alpha=4\mu^{-1}(\mu/4\pm\delta^2\mu^{-1}\mp c^2)h,\ \
     \beta=4\mu^{-\frac{3}{2}}\sqrt{\mu/4\pm\delta^2\mu^{-1}\mp c^2}\
     c_0.
  \ee

The proof of Theorem \ref{th3} (ii) is similar to that of Theorem
\ref{th2} (ii). Let $x(t)$ ($0\le t\le 2\pi/\sqrt{\mu}$) be an
arbitrary closed geodesic of $h$ connecting the two points $P,Q$
on which the gradient  $\nabla c$ vanishes. Along $x(t)$, we have
$$c=\mu^{-\frac{1}{2}}\delta \cos(\sqrt{\mu}\ t),\ \
h(x(t),x'(t))=1,\ \ c_0(x(t),x'(t))=c'(t).$$
 Since $\beta$ is closed, it follows from (\ref{y45}) that
 the $F$-length $L_1$ of the geodesic $x(t)$ is given by
  \beqn
  L_1&=&\int_0^{\frac{2\pi}{\sqrt{\mu}}} F_{\epsilon}^{\pm}(x(t),x'(t))dt
  =\int_0^{\frac{2\pi}{\sqrt{\mu}}} F_0^{\pm}(x(t),x'(t))dt\\
  &=&\int_0^{\frac{2\pi}{\sqrt{\mu}}}\Big\{
  \alpha(x(t),x'(t))\pm\frac{\beta^2(x(t),x'(t))}{\alpha(x(t),x'(t))}\Big\}dt\\
  &=&\int_0^{\frac{2\pi}{\sqrt{\mu}}}
  \Big\{1\pm\frac{4\delta^2}{\mu^2}\mp\frac{4\delta^2}{\mu^2}\cos(2\sqrt{\mu}\
  t)\Big\}dt\\
 &=&\frac{2\pi}{\sqrt{\mu}}\pm
  \frac{8\pi}{\mu^2\sqrt{\mu}}\ \delta^2.
  \eeqn
Thus we obtain (\ref{y12}) for the $F$-length of the family of
closed geodesics collecting  $P,Q$.

Now let $z(t)$ ($0\le t\le 2\pi/\sqrt{\mu}$) be an arbitrary
closed geodesic of $(S^n,h)$ lying on the hypersurface $c=0$.  By
$c(z(t))=0$, $c_0(z(t),z'(t))=0$ and $h(z(t),z'(t))=1$, it is
easily seen that the $F$-length $L_2$ of the geodesic $z(t)$ is
given by
 \beqn
 L_2&=&\int_0^{\frac{2\pi}{\sqrt{\mu}}}
 F_{\epsilon}^{\pm}(x(t),x'(t))dt=\int_0^{\frac{2\pi}{\sqrt{\mu}}}\alpha(z(t),z'(t))dt\\
 &=&\int_0^{\frac{2\pi}{\sqrt{\mu}}}\big(1\pm\frac{4\delta^2}{\mu^2}\big)dt=\frac{2\pi}{\sqrt{\mu}}\pm
  \frac{8\pi}{\mu^2\sqrt{\mu}}\ \delta^2.
 \eeqn
This gives (\ref{y12}) for the $F$-length of the closed geodesics
lying on the hypersurface $c=0$.

Finally, we come to the proof of item (iii) of this theorem. Plug
$k_1=\pm 2,k_2=0,k_3=\mp 3$, (\ref{y43}) and (\ref{y44}) into the
scalar flag curvature $K_{\epsilon}^{\pm}=R(\beta/\alpha,c^2)$ for
$F=F_{\epsilon}^{\pm}$, where $R=R(s,t)$ is given by (\ref{y9}),
and then we obtain
 \be\label{y46}
   K_{\epsilon}^{\pm}=\frac{6\mu(\epsilon^2\mp4)(1\mp
s^2)^2c^2+(\mu^2\pm4 \delta^2)(\pm \epsilon s^3\pm 6 s^2+3\epsilon
s+2)(1+\epsilon s\pm s^2)}{128\mu^{-2}(\mu/4\pm\delta^2\mu^{-1}\mp
c^2)^3(1+\epsilon s\pm s^2)^4}.
 \ee

\noindent {\bf Case (A) :}\
 For $F=F_{2}^{+}$, by (\ref{y46}) we
get
 $$
 K_2^+
=\frac{\xi\mu^3}{16}\big[(1+s)(\xi-c^2)\big]^{-3},\ \
\big(\xi:=\delta^2\mu^{-1}+\mu/4\big).
$$
By Theorem \ref{th2} (iii), there is a function $R=R(s,t)$ of two
variables such that $K_2^+=R(\beta/\alpha,c^2)$. By (\ref{y42}) in
Remark \ref{rem42}, $R$ can be written as
$\widetilde{R}=\widetilde{R}(s,t)$ defined on $\widetilde{D}$ with
 \beqn
 &&\ \ \ \ \widetilde{R}=\frac{\xi\mu^3}{16}\big[(1+s)\ \frac{\mu}{4(1-t^2)}\big]^{-3}
 =4\xi\big(\frac{1+s}{1-t^2}\big)^{-3},\\
 &&
 \widetilde{D}=\{(s,t)|0\le t\le t_o,\ \ |s|\le t\},\ \
 \big(t_o:=\frac{2\delta}{\sqrt{4\delta^2+\mu^2}}\big).
 \eeqn
It is clear that
 $$Min_{(s,t)\in
 \widetilde{D}}\Big(\frac{1+s}{1-t^2}\Big)=Min_{t\in [0,
 t_o]}\Big(\frac{1}{1+t}\Big)=\frac{1}{1+t_o},$$
 $$Max_{(s,t)\in
 \widetilde{D}}\Big(\frac{1+s}{1-t^2}\Big)=Max_{t\in [0,
 t_o]}\Big(\frac{1}{1-t}\Big)=\frac{1}{1-t_o}.
$$
Now we can easily obtain the maximal and minimal values of $K_2^+$
on $S^n$ given by (\ref{y13}).

\

\noindent {\bf Case (B) :}  \ For $F=F_{0}^{+}$, by (\ref{y46}) we
get
 $$
 K_0^+
=\frac{\mu^5\big[-12\mu(1-s^2)^2c^2+(\mu^2+4\delta^2)(1+s^2)(1+3s^2)\big]}{(\mu^2+4\delta^2-4\mu
c^2)^3(1+s^2)^4}.
$$
Then in Remark \ref{rem42}, $\widetilde{R}=\widetilde{R}(s,t)$ and
$\widetilde{D}$ are given by
 \beq
&&\widetilde{R}=\frac{(1-t^2)^2\big[2(\mu^2+4\delta^2)(1-5s^2)t^2+3\mu^2s^4+4(\mu^2+10\delta^2)s^2+\mu^2-8\delta^2\big]}{\mu(1+s^2)^4},\label{y47}\\
 &&\hspace{1.5cm}\widetilde{D}=\{(s,t)|0\le t\le t_o,\ \ |s|\le t\},\ \
 \big(t_o:=\frac{2\delta}{\sqrt{4\delta^2+\mu^2}}\big).\label{y48}
 \eeq

(a1). For the function $\widetilde{R}$ defined by (\ref{y47}), a
direct computation shows that $d\widetilde{R}\ne 0$ for $(s,t)$
belonging to the interior of $\widetilde{D}$ defined by
(\ref{y48}). So the maximal and minimal values of $\widetilde{R}$
are taken on the boundary of $\widetilde{D}$.

(a2). On the boundary $|s|=t$, we have
 $$
 \varphi(t):=\widetilde{R}(\pm
 t,t)=\frac{(1-t^2)^3\big[(7\mu^2+40\delta^2)t^2+\mu^2-8\delta^2\big]}{\mu(1+t^2)^4},
 \ \ 0\le t\le t_o.
 $$
It is easy to see that $\varphi'(t)=0$ has a unique solution
$t_1\in(0,t_o)$, where
 $$
 t_1=\frac{\sqrt{2}\ \delta}{\sqrt{\mu^2+6\delta^2}}.
 $$

 (a3). On the boundary $t=t_o$, we have
  $$
 \psi(s):=\widetilde{R}(s,t_o)=\frac{\mu^5}{2(\mu^2+4\delta^2)^2}\frac{1+3s^2}{(1+s^2)^3},
 \ \ \ s\in(-t_o,t_o).
  $$
It is easy to see that $\psi'(s)=0$ has a unique solution $s_1$,
where $s_1=0$.

Summing up on (a1)--(a3), the maximal and minimal values of
$\widetilde{R}(s,t)$ on $\widetilde{D}$ are obtained from the
following four values
 $$\widetilde{R}(0,0),\ \ \widetilde{R}(\pm t_1,t_1),\ \ \widetilde{R}(\pm
 t_o,t_o),\ \ \widetilde{R}(s_1,t_o).$$
Therefore, the minimal and maximal values of the scalar flag
curvature $K_0^+$ on $S^n$ are given respectively by
 $$
 Min(K_0^+)=\widetilde{R}(0,0)=\frac{\mu^2-8\delta^2}{\mu},\ \ \
 Max(K_0^+)=\widetilde{R}(\pm t_1,t_1)=\frac{(\mu^2+4\delta^2)^4}{\mu(\mu^2+8\delta^2)^3}
 $$
 Then we obtain the proof of (\ref{y14}).

 \

\noindent {\bf Case (C) :}\  For $F=F_{0}^{-}$, first note that
$\mu^2>12\delta^2$ in the following discussion (see Remark
\ref{rem51} below). In this case, by (\ref{y46}) we get
 $$
 K_0^-
=\frac{\mu^5\big[12\mu(1+s^2)^2c^2+(\mu^2-4\delta^2)(1-s^2)(1-3s^2)\big]}{(\mu^2-4\delta^2+4\mu
c^2)^3(1-s^2)^4}.
$$
Then in Remark \ref{rem42}, $\widetilde{R}=\widetilde{R}(s,t)$ and
$\widetilde{D}$ are given by
 \beq
&&\widetilde{R}=\frac{(1+t^2)^2\big[2(4\delta^2-\mu^2)(1+5s^2)t^2+3\mu^2s^4-4(\mu^2-10\delta^2)s^2+\mu^2+8\delta^2\big]}{\mu(1-s^2)^4},\label{y49}\\
 &&\hspace{1.5cm}\widetilde{D}=\{(s,t)|0\le t\le t_o,\ \ |s|\le t\},\ \
 \big(t_o:=\frac{2\delta}{\sqrt{\mu^2-4\delta^2}}\big).\label{y50}
 \eeq

(b1). For the function $\widetilde{R}$ defined by (\ref{y49}), a
direct computation shows that $d\widetilde{R}\ne 0$ for $(s,t)$
belonging to the interior of $\widetilde{D}$ defined by
(\ref{y50}). So the maximal and minimal values of $\widetilde{R}$
are taken on the boundary of $\widetilde{D}$.

(b2). On the boundary $|s|=t$, we have
 $$
 \lambda(t):=\widetilde{R}(\pm
 t,t)=\frac{(1+t^2)^3\big[(40\delta^2-7\mu^2)t^2+\mu^2+8\delta^2\big]}{\mu(1-t^2)^4},
 \ \ 0\le t\le t_o.
 $$
It is easy to see that $\lambda'(t)=0$ has a unique solution
$t_1\in(0,t_o)$, where
 $$
 t_1=\frac{\sqrt{2}\ \delta}{\sqrt{\mu^2-6\delta^2}}.
 $$

 (b3). On the boundary $t=t_o$, we have
  $$
 \chi(s):=\widetilde{R}(s,t_o)=\frac{\mu^5}{2(\mu^2-4\delta^2)^2}\frac{1-3s^2}{(1-s^2)^3},
  \ \ \ s\in(-t_o,t_o).
  $$
It is easy to see that $\chi'(s)=0$ has a unique solution
$s_1\in(-t_o,t_o)$, where $s_1=0$.

Summing up on (b1)--(b3), the maximal and minimal values of
$\widetilde{R}(s,t)$ on $\widetilde{D}$ are obtained from the
following four values
 $$\widetilde{R}(0,0),\ \ \widetilde{R}(\pm t_1,t_1),\ \ \widetilde{R}(\pm
 t_o,t_o),\ \ \widetilde{R}(s_1,t_o).$$
Therefore, the minimal and maximal values of the scalar flag
curvature $K_0^-$ on $S^n$ are given respectively by
 $$
 Min(K_0^-)=\widetilde{R}(\pm
 t_o,t_o)=\frac{\mu^5(\mu^2-16\delta^2)}{(\mu^2-8\delta^2)^3},\ \ \
 Max(K_0^-)=\widetilde{R}(\pm
 t_1,t_1)=\frac{(\mu^2-4\delta^2)^4}{\mu(\mu^2-8\delta^2)^3}.
 $$
 Then we obtain the proof of (\ref{y15}).  \qed

\begin{rem}\label{rem51}
In Theorem \ref{th2},  it is easily seen from (\ref{y44}) that
$F=F_{2}^+$ and $F=F_{0}^+$ are
 regular if and only if the constants $\mu,\delta$ satisfy $\mu>0$ and
 $\delta\ge 0$; $F=F_{0}^-$ is regular if and only if the constants
  $\mu,\delta$ satisfy $\mu^2>12\delta^2$. The former follows from
(\ref{y44}) and $b^2<1$, and the latter holds since by (\ref{y44})
we get $b^2<1/2\Leftrightarrow c^2>(12\delta^2-\mu^2)/(12\mu)$,
and thus we have $12\delta^2-\mu^2<0$ from $Min_{x\in S^n}c(x)=0$.
\end{rem}
 ,

\vspace{0.5cm}

\noindent Guojun Yang \\
Department of Mathematics \\
Sichuan University \\
Chengdu 610064, P. R. China \\
 yangguojun@scu.edu.cn

\end{document}